\documentclass[12pt,reqno,twoside,a4paper]{article}
\usepackage{longtable}
\usepackage{amsmath,amssymb,amsfonts,graphicx}
\usepackage[mathscr]{eucal}
\usepackage[T2A]{fontenc}
\usepackage[cp866]{inputenc}
\usepackage[english]{babel}
\def\{{\protect\lbrace}
\def\}{\protect\rbrace}
\def\S{$\mathchar"278$}

\newcommand{\End}{\operatorname{End}}
\newcommand{\Hom}{\operatorname{Hom}}

\begin{document}
\begin{center}

\textbf{Automorphism-Liftable Modules}
\end{center}

\hfill {\sf A.A. Tuganbaev}

\hfill National Research University "MPEI"

\hfill Lomonosov Moscow State University

\hfill e-mail: tuganbaev@gmail.com

\textbf{Abstract.} In this paper, we describe all automorphism-liftable torsion modules over non-primitive hereditary Noetherian prime rings. We also study automorphism-liftable non-torsion modules over not necessarily commutative Dedekind prime rings.

The work is supported by Russian Scientific Foundation, project 16-11-10013.


2000 MATHEMATICS SUBJECT CLASSIFICATION 16D40; 16D80; 16N80

\section{Introduction}

All considered rings are associative and contain the non-zero identity element. Writing expressions of the form ``$A$ is a non-primitive ring or a Noetherian ring'' we mean that $A$ is not a right and left primitive ring or the both $A_A$ and $_AA$ are Noetherian.

\textbf{1.1. Automorphism-liftable, strongly automorphism-liftable, endomorphism-liftable, and quasi-projective modules.}

A module $M$ is said to be \textsf{automorphism-liftable}\footnote{The notion of an automorphism-liftable module is dual to the notion of an automorphism-extendable module studued in \cite{Tug13} and \cite{Tug16}. A module $M$ is said to be \textsf{automorphism-extendable} if every automorphism of any its submodule can be extended to an endomorphism of the module $M$.} (resp., \textsf{strongly automorphism-liftable}) if for any epimorphism $h\colon M\to \bar{M}$ and every automorphism $\bar{f}$ of the module $\bar{M}$, there exists an endomorphism (resp., automorphism) $f$ of the module $M$ with $\bar{f}h=hf$.

A module $M$ is said to be \textsf{endomorphism-liftable} if for any epimorphism $h\colon M\to \bar{M}$ and every endomorphism $\bar{f}$ of the module $\bar{M}$, there exists an endomorphism $f$ of the module $M$ with $\bar{f}h=hf$. In the above definitions, without loss of generality, we can assume that $\bar{M}$ is an arbitrary factor module of the module $M$ and $h\colon M\to \bar{M}$ is the natural epimorphism.

It is clear that all strongly automorphism-liftable modules and all endomorphism-liftables modules are automorphism-liftable. 

Endomorphism-liftable modules and Abelian groups were studied in many papers under various names; e.g., see \cite{Jan73}, \cite{Mis72}, \cite{Tug78}, \cite{Tug80}, \cite{Tug80b}, \cite{Tug16}. In particular, endomorphism-liftable Abelian groups were studied in \cite{Mis72} and \cite{Jan73}; endomorphism-liftable modules over non-primitive hereditary Noetherian prime rings were studied in \cite{Tug80}, \cite{Tug80b}. 

\textbf{1.2. Automorphism-liftable $\mathbb{Z}$-modules.}

Any quasi-cyclic Abelian group $\mathbb{Z}(p^{\infty})$ is an endomorphism-liftable (automorphism-liftable) non-quasi-projective $\mathbb{Z}$-module.

The ring of integers $\mathbb{Z}$ is an automorphism-liftable $\mathbb{Z}$-module which is not strongly automorphism-liftable. Indeed, let $\bar{f}$ be an automorphism of the simple $\mathbb{Z}$-module $\mathbb{Z}/5\mathbb{Z}$ such that $\bar{f}$ multiplies all elements of this module by $3$. Since the only non-identity automorphism of the module $\mathbb{Z}_{\mathbb{Z}}$ coincides with the multiplication by $-1$, the projective module $\mathbb{Z}_{\mathbb{Z}}$ is not strongly automorphism-liftable.

A.P.Mishina \cite{Mis72} completely described strongly automorphism-liftable Abelian groups, i.e., strongly automorphism-liftable $\mathbb{Z}$-modules. It follows from this description that strongly automorphism-liftable $\mathbb{Z}$-modules are torsion automorphism-liftable $\mathbb{Z}$-modules.
In \cite{AbyQ18}, A.N.Abyzov and T.C.Quynh described torsion automorphism-liftable\footnote{In this paper, automorphism-liftable modules are called \textsf{dually automorphism-extendable}. Automorphism-liftable modules are also studied in \cite{SelS18}.} $\mathbb{Z}$-modules; it follows from this description and results of A.P.Mishina \cite{Mis72} that the strongly automorphism-liftable $\mathbb{Z}$-modules coincide with torsion automorphism-liftable $\mathbb{Z}$-modules.

\textbf{1.3. Non-primitive hereditary Noetherian prime rings.} The ring $\mathbb{Z}$ is a very partial case of a non-primitive hereditary Noetherian prime ring.\footnote{A module is said to be \textsf{hereditary} if all its submodules are projective.} A ring $A$ is said to be \textsf{right bounded} (resp., \textsf{left} bounded) if every its essential right (resp., left) ideal contains a non-zero ideal of the ring $A$. If $A$ is a non-primitive hereditary Noetherian prime ring, then $A$ is not a right or left Artinian; see \cite{Len73}. Every hereditary Noetherian prime ring $A$ is a primitive ring or a bounded ring and if $A$ is a primitive bounded ring, then $A$ is a simple Artinian ring; see \cite{Len73}.

Let $A$ be a Noetherian prime ring. It is well known that the ring $A$ has the simple Artinian classical ring of fractions $Q$. An ideal $B$ of the ring $A$ is called an \textsf{invertible} ideal if there exists a subbimodule $B^{-1}$ of the bimodule ${}_AQ_A$ such that $BB^{-1}=B^{-1}B=A$. 
The maximal elements of the set of all proper invertible ideals of the ring $A$ are called \textsf{maximal invertible ideals} of the ring $A$. The set of all maximal invertible ideals of the ring $A$ is denoted by $\mathbb{P}(A)$. If $P\in \mathbb{P}(A)$, then the submodule $\{m\in M\;\vert \;mP^n=0$, $n=1,2, \ldots \}$ is called the \textsf{$P$-primary component} of the module $M$; it is denoted by $M(P)$. If $M=M(P)$ for some $P\in \mathbb{P}(A)$, then $M$ is called a \textsf{primary} module or a \textsf{$P$-primary} module.

In connection with 1.2, we prove Theorem 1.4 which is the first main result of this paper. This theorem generalizes the description of torsion automorphism-liftable $\mathbb{Z}$-modules from \cite{AbyQ18} for the case of singular modules over non-primitive hereditary Noetherian prime rings.\footnote{Over any Noetherian prime ring $A$, the singular modules coincide with the torsion modules, where a module $M_A$ is said to be \textsf{singular} (resp., \textsf{torsion} if the annihilator of every its element is an essential right ideal of the ring $A$ (resp., contains a non-zero-divisor of the ring $A$).} 

\textbf{1.4. Theorem.} If $A$ is a non-primitive hereditary Noetherian prime ring and $M$ is a singular right $A$-module, then $M$ is automorphism-liftable if and only if every $P$-primary component\footnote{Necessary definitions are given at the end of Introduction.} $M(P)$ of the module $M$ is either a projective $A/r(M(P))$-module or a uniserial injective module $M(P)$ such that all proper submodules are cyclic and form a countable chain
$$
0=x_0R\subsetneq x_1R\subsetneq \ldots , 
$$
all subsequent factors of this chain are simple modules and there exists a positive integer $n$ such that $x_kA\cong x_{k+n}/x^nA$ and $M(P)/x_kA\cong M(P)/x_{k+n}A$ for all $k=0,1,2,\ldots$.

\textbf{1.5. Remark.} In \cite{Tug80b}, arbitrary endomorphism-liftable modules over non-primitive hereditary Noetherian prime rings are described.

\textbf{1.6. Non-primitive Dedekind prime rings.} A hereditary Noetherian prime ring $A$ is  called a \textsf{Dedekind prime} ring (see \cite[\S 5.2]{MccR87}) if every its non-zero ideal is invertible in the simple Artinian ring of fractions of the ring $A$. Any hereditary Noetherian prime PI ring is a bounded ring. In particular, all commutative Dedekind domains (e.g., the ring $\mathbb{Z}$) and matrix rings over commutative Dedekind domains are non-primitive Dedekind prime rings. Other examples of Dedekind prime rings are given in \cite[\S~5.2,~\S~5.3]{MccR87}.

The second main result of this paper is Theorem 1.7, where the description of singular automorphism-liftable modules from Theorem 1.4 is specified in the case, where $A$ is a non-primitive Dedekind prime ring.

\textbf{1.7. Theorem.} If $A$ is a non-primitive Dedekind prime ring and $M$ is a singular right $A$-module, then the module $M$ is automorphism-liftable if and only if every $P$-primary component $M(P)$ of the module $M$ is either the direct sum of isomorphic cyclic uniserial modules of finite length or a uniserial injective module $M(P)$ such that all proper submodules are cyclic and form a countable chain
$$
0=x_0A\subsetneq x_1A\subsetneq \ldots 
$$
such that all subsequent factors of this chain are isomorphic simple modules and $M(P)/x_kA\cong M(P)/x_{k+1}A$ for all $k=0,1,2,\ldots$.

\textbf{1.8. Theorem \cite{Sin74}.} If $R$ is a non-primitive Dedekind prime ring with simple Artinian ring of fractions $Q$, then the following conditions are equivalent.

\textbf{1)} $Q_R$ is a quasi-projective module.

\textbf{2)} $R=D_n$, where $D$ is a local principal right (left) domain which is complete in the topology defined by powers of the Jacobson radical $J(D)$.

In this case, $R$ is a Dedekind prime ring, $J(R)$ is a maximal ideal, and $_RQ$ is a quasi-projective module.

A ring $R$ is called a \textsf{special Dedekind prime} ring if $R$ satisfies Theorem 1.8.

The third main result of this paper is Theorem 1.9, where all automorphism-liftable $A$-modules are described in the case, where $A$ is a non-primitive Dedekind prime ring with $\displaystyle{\frac{1}{2}}$.

\textbf{1.9. Theorem.} Let $A$ be a non-primitive Dedekind prime ring with $\frac{1}{2}\in A$. A right $A$-module $M$ is automorphism-liftable if and only if one of the following three conditions holds:

\textbf{a)} If $M$ is a torsion module, then $M$ is endomorphism-liftable if and only if every primary component of the module $M$ is either an indecomposable injective module or a projective $A/r(M)$-module.

\textbf{b)} If $M$ is a mixed module, then $M$ is endomorphism-liftable if and only if $M=T\oplus F$, where $T$ is a torsion injective module such that all primary components are indecomposable and $F$ is a finitely generated projective module.

\textbf{c)} If $M$ is a torsion-free module, then $M$ is endomorphism-liftable if and only if either $M$ is projective or $A$ is a special ring with classical ring of fractions $Q$ and $M=E\oplus F$, where $E$ is a minimal right ideal of the ring $Q$, $F$ is a finitely generated projective module, and $n$ is a positive integer.

We give some necessary definitions and notation.

Let $A$ be a ring and $M$ a right $A$-module.

We denote by $J(A)$ the Jacobson radical of $A$. We denote by $r(X)$ the annihilator in the ring $A$ of the subset $X$ of the module $M$. 

We denote by $T(M)$ the set of all elements of $M$ whose annihilators contain a non-zero-divisor; the set $T(M)$ is called the \textsf{torsion part} of the module $M$. If $T(M)=0$ (resp., $0\ne T(M)\ne M$), then $M$ is said to be \textsf{torsion-free module} (resp., \textsf{mixed}). Every module is a torsion module, or a torsion-free module, or a mixed module.

A module $M$ is said to be \textsf{projective with respect to} the module $N$ (or \textsf{$N$-projective}) if for any epimorphism $h\colon N\to \bar N$ and every homomorphism $\bar f\colon M\to \bar N$, there exists a homomorphism $f\colon M\to N$ with $hf=\bar f$. Thus, quasi-projective modules coincide with the modules which are projective with respect to itself. Clearly, in the definition of the relative projectivity, it is sufficient to consider only the case, where $h\colon N\to \bar N$ is the natural epimorphism of the module $N$ onto its arbitrary factor module $\bar N$.

A module is said to be \textsf{uniserial} if any two of its submodules are comparable with respect to inclusion. For a module $M$, a submodule $X$ of $M$ is said to be \textsf{essential} (in $M$) if $X$ has the non-zero intersection with any non-zero submodule of the module $M$. A module $M$ is said to be \textsf{finite-dimensional} if $M$ does not contain an infinite direct sum of non-zero submodules.

\section{Proof of Theorem 1.4 and Theorem 1.7} \label{section2}

\textbf{2.1. Lemma.} If $A$ is a ring and $M$ is a right $A$-module, then $M$ is an automorphism-liftable (resp., idempotent-liftable, quasi-projective,) $A$-module if and only if $M$ is an automorphism-liftable (resp., idempotent-liftable; quasi-projective) $A/r(M)$-module.

Lemma 2.1 is directly verified.

\textbf{2.2. Lemma.} Let $M=\oplus _{i\in I}M_i$ be a module and $N= \oplus _{i\in I}(N\bigcap M_i)$ for any submodule $N$ of the module $M$. Then:

\textbf{1)} $\Hom (M_i,M_j)=0$ for any $i\ne j$ in $I$ (therefore, all the $M_i$ are fully invariant submodules in $M$);

\textbf{2)} if $N$, $P$ and $Q$ are three submodules of the module $M$, then the relation $N=P+Q$ is equivalent to the property that $N\bigcap M_i=P\bigcap M_i+Q\bigcap M_i$ for all $i\in I$;

\textbf{3)} the module $M$ is automorphism-liftable if and only if each of the modules $M_i$ are automorphism-liftable.

\textbf{Proof.} \textbf{1) and 2).} The assertions are proved in \cite[Lemma 2.1(1),(2)]{Tug01}.

\textbf{3.} The assertion follows from \textbf{1)} and \textbf{2)}.~\hfill$\square$

\textbf{2.3. Lemma \cite[Proposition 6, Lemma 2]{AbyQ18}.}
If an automorphism-liftable module $M$ is the direct sum of some modules $X$ and $Y$, then $X$, $Y$ are automorphism-liftable modules which are projective with respect to each other.

\textbf{2.4. Lemma.} If $A$ is a non-primitive hereditary Noetherian prime ring and $M$ is a primary right $A$-module, then the following conditions are equivalent.

\textbf{1)} $M$ is an automorphism-liftable module;

\textbf{2)} $M$ is an endomorphism-liftable module;

\textbf{3)} for any direct decomposition $M=M_1\oplus M_2$, the module $M_1$ is projective with respect to the module $M_2$;

\textbf{4)} $M$ is a projective $A/r(M)$-module or an indecomposable injective $A$-module.

\textbf{Proof.} The equivalence of conditions 2), 3) and 4) is proved in \cite[Lemma 13]{Tug80b}. 

2)\,$\Rightarrow$\,1). The assertion is always true.

1)\,$\Rightarrow$\,2). The assertion follows from Lemma 2.4.~\hfill$\square$

\textbf{2.5. Theorem.} Let $A$ be a non-primitive hereditary Noetherian prime ring, $M$ a torsion right $A$-module, $\{M_i\}$ the set of all primary components of the module $M$. Then the module $M$ is automorphism-liftable if and only if every primary component $M_i$ of the module $M$ is a projective $A/r(M)$-module or an indecomposable injective $A$-module.

\textbf{Proof.} For an arbitrary torsion $A$-module $M$ and every its submodule $N$, we have $N= \oplus _{i\in I}(N\bigcap M_i)$; see \cite[Lemma 2.2(1)]{Tug01}. Therefore, our assertion follows from Lemma 2.4.~\hfill$\square$

\textbf{2.6. Remark.} If $A$ is a non-primitive hereditary Noetherian prime ring, then the structure of indecomposable injective torsion $A$-modules is known; e.g., see \cite{Sin74} and \cite{Sin75}. Namely, the indecomposable injective torsion $A$-modules coincide with the primary modules $M$ such that all proper submodules of $M$ are cyclic uniserial primary modules and form a countable chain
$$
0=x_0R\subsetneq x_1R\subsetneq \ldots ,
$$
all subsequent factors of this chain are simple modules and there exists a positive integer $n$ such that $x_kA\cong x_{k+n}/x^nA$ and $M(P)/x_kA\cong M(P)/x_{k+n}A$ for all $k=0,1,2,\ldots$.\\
If $A$ is a non-primitive Dedekind prime ring, then $n=1$.

\textbf{2.7. Remark.} If $A$ is a non-primitive Dedekind prime ring and $M$ is a torsion $A$-module, then $M$ is quasi-projective if and only if every primary component of the module $M$ is a direct sum of isomorphic cyclic modules of finite length. \cite[Theorem 15]{Sin75}

\textbf{2.8. Completion of the proof of Theorem 1.4 and Theorem 1.7.}\\
Theorem 1.4 and Theorem 1.7 follow from Theorem 2.5, Remark 2.6 and Remark 2.7.~\hfill$\square$

\section{Proof of Theorem 1.9} \label{section3}

\textbf{3.1. Idempotent-liftable modules and $\pi$-projective modules.}\\
A module $M$ is said to be \textit{idempotent-liftable} if for any epimorphism $h\colon M\to \bar{M}$ and every idempotent endomorphism $\bar{f}$ of the module $\bar{M}$, there exists an endomorphism $f$ of the module $M$ with $\bar{f}h=hf$. Without loss of generality, we can assume that $\bar{M}$ is an arbitrary factor module of the module $M$ and $h\colon M\to \bar{M}$ is the natural epimorphism. 

A module $M$ is said to be \textsf{$\pi$-projective} if for any its submodules $X$ and $Y$ with $X+Y=M$, there exist endomorphisms $f$ and $g$ of the module $M$ with $f+g=1_M$, $f(M)\subseteq X$ and $g(M)\subseteq Y$ (see \cite[p.359]{Wis91}).\\
The idempotent-liftable modules coincide with the $\pi$-projective modules; see Lemma 3.2.

It is clear that every endomorphism-liftable module is idempotent-liftable. If $A$ is a uniserial principal right (left) ideal domain with division ring of fractions $Q$ which is not complete with respect to the $J(A)$-adic topology, then $Q$ is an idempotent-liftable $A$-module which is not endomorphism-liftable.

The class of all idempotent-liftable modules contains all modules such that all their factor modules are indecomposable, all projective modules, all uniserial modules, and all local modules. In particular, all cyclic modules over local rings and all free modules are idempotent-liftable. 

\textbf{3.2. Lemma.} For a module $M$, the following conditions are equivalent.

\textbf{1)} $M$ is an idempotent-liftable module.

\textbf{2)} $M$ is a $\pi$-projective module. 

\textbf{Proof.} 1)\,$\Rightarrow$\,2). Let $M=X+Y$ be an idempotent-liftable module, $N=X\cap Y$, $\bar M=M/N$, $\bar X=X/N$, $\bar f$, $\bar g_1$ natural projections of the module $\bar M=\bar X\oplus \bar Y$ onto the components $\bar X$, $\bar Y$, respectively, and $h$ the natural epimorphism from the module $M$ onto $\bar M$. Since $M$ is an idempotent-liftable module, there are two endomorphisms $f$, $g_1$ of the module $M$ such that $\bar f h=h f$, $\bar g_1 h=h g_1$. It is easy to see that $fM\subseteq X$, $g_1M\subseteq Y$. Since $f+g_1$ coincides  with the identity mapping on $\bar M$, we have $(1_M-f-g_1)M\subseteq N$. We set $g=l_M-f$. Then $fM\subseteq X$, $gM\subseteq (1-f-g_1)M + g_1M\subseteq N+Y=Y$. 

2)\,$\Rightarrow$\,1). Let $\bar{M}=M/N$ be an arbitrary factor module of the module $M$, $h\colon M\to \bar{M}$ the natural epimorphism, $\bar{f}$ an idempotent endomorphism of the module $\bar{M}$, $\bar{X}=X/N=\bar{f}(\bar{M})$, $\bar{Y}=Y/N=(1-\bar{f})(\bar{M})$, where $X$, $Y$ are complete pre-images of the modules $\bar X$ and $\bar Y$ in $M$, respectively. Then $\bar M=\bar X\oplus \bar Y$, $M=X+Y$ and $N=X\cap Y$. Since $M$ is a $\pi$-projective module, there exist homomorphisms $f\colon M\to X$ and $g\colon M\to Y$ with $f+g=1_M$. Then $\bar{f} h=hf$ and $M$ is an idempotent-liftable module.~\hfill$\square$

\textbf{3.3. Lemma.} If $A$ is a ring with $\displaystyle{\frac{1}{2}}$ and $M$ is a right $A$-module, then every idempotent endomorphism of the module $M$ is the sum of two automorphisms of the module $M$; in particular, every automorphism-liftable (resp., automorphism-extendable\footnote{A module $M$ is said to be \textsf{automorphism-extendable} if every automorphism of any its submodule can be extended to an endomorphism of the module $M$.}) right $A$-module is a idempotent-liftable (resp., idempotent-extendable\footnote{A module $M$ is said to be \textsf{idempotent-extendable} if every idempotent endomorphism of any its submodule can be extended to an endomorphism of the module $M$.}).

\textbf{Proof.} Let $f$ be an idempotent endomorphism of the module $M$. Then $M=X\oplus Y$, where $X=f(M)$ and $Y=(1_M-f)(M)$. We denote by $u$ an automorphism of the module $M$ such that $u(x+y)=x-y$ for $x\in X$, $y\in Y$. Then $f=\frac{1}{2}\cdot 1_M+\frac{1}{2}\cdot u$, where $\frac{1}{2}\cdot 1_M$ and $\frac{1}{2}\cdot u$ are automorphisms of the module $M$.~\hfill$\square$

\textbf{3.4. Theorem \cite[Theorem 1]{Tug01}.} A module $M$ over a non-primitive Dedekind prime ring $A$ is $\pi$-projective if and only if one of the following three conditions holds:

\textbf{a)} $M$ is a torsion module such that every primary
component is either an indecomposable injective module or
direct sum of isomorphic cyclic modules of finite length;

\textbf{b)} $M=T\oplus F$, where $T$ is a non-zero injective torsion module such that every primary component is an indecomposable module and $F$ is a non-zero finitely generated projective module;

\textbf{c)} $M$ is a projective module or there exist two positive integers $k$ and $n$ such that the ring $A$ is isomorphic to the ring of all $k\times k$ matrices $D_k$ over some uniserial principal right (left) ideal domain $D$, $M=X\oplus Y$, $X$ is the finite direct sum of non-zero injective indecomposable torsion-free modules $X_1, \ldots ,X_n$,~ $Y$ is a finitely generated projective module, and either $n=1$ or $n\ge 2$ and $D$ is a complete domain.

\textbf{3.5. Theorem \cite{Tug80b}.} Let $A$ be a non-primitive hereditary Noetherian prime ring and $M$ a right $A$-module.

\textbf{a)} If $M$ is a torsion module, then $M$ is endomorphism-liftable if and only if every primary component of the module $M$ is either an indecomposable injective module or a projective $A/r(M)$-module.

\textbf{b)} If $M$ is a mixed module, then $M$ is endomorphism-liftable if and only if $M=T\oplus F$, where $T$ is a torsion injective module such that all primary components are indecomposable and $F$ is a finitely generated projective module.

\textbf{c)} If $M$ is a torsion-free module, then $M$ is endomorphism-liftable if and only if either $M$ is projective or $A$ is a special ring with classical ring of fractions $Q$ and $M=E^n\oplus F$, where $E$ is a minimal right ideal of the ring $Q$,~ $F$ is a finitely generated projective module, and $n$ is a positive integer.

\textbf{3.6. Theorem \cite[Theorem 2]{Tug80}.} If $A$ is a non-primitive hereditary Noetherian prime ring, then a right $A$-module $M$ is quasi-projective if and only if either $M$ is a torsion module and every its primary component is a projective $A/r(M)$-module, or $M$ is projective, or $A$ is a special ring and $M=E\oplus F$, where $E$ is an injective finite-dimensional torsion-free module and $F$ is a finitely generated projective module.

\textbf{3.7. Lemma.} If $M$ is an automorphism-liftable module with local endomorphism ring $\End M$, then $M$ is an endomorphism-liftable module.

\textbf{Proof.} Let $h\colon M\to \bar{M}$ be an epimorphism and $\bar{f}$ an endomorphism of the module $\bar{M}$. If $\bar{f}$ is an automorphism of the module $\bar{M}$, then it follows from the assumption that there exists an endomorphism $f$ of the module $M$ that $\bar{f}h=hf$.

We assume that $\bar{f}$ is not an automorphism of the module $M$. By assumption, the ring $\End M$ is local. Therefore, $1_{\bar{M}}-\bar{f}$ is an automorphism of the module $\bar{M}$. Since $M$ is an automorphism-liftable module, there exists an endomorphism $g$ of the module $M$ such that $(1_{\bar{M}}-\bar{f})h=hg$. We denote by $f$ the endomorphism $1-g$ of the module $M$. Since $h-1_{\bar{M}}h=0$, we have
$$
hf=h-hg=h-1_{\bar{M}}h+\bar{f}h=\bar{f}h.
$$
Therefore, $M$ is an endomorphism-liftable module.~\hfill$\square$

\textbf{3.8. Theorem.} A module $M$ over a non-primitive Dedekind prime ring $A$ is an automorphism-liftable, idempotent-liftable module if and only if one of the following three conditions holds:

\textbf{a)} If $M$ is a torsion module, then $M$ is a endomorphism-liftable if and only if every primary component of the module $M$ is either an indecomposable injective module or a projective $A/r(M)$-module.

\textbf{b)} If $M$ is a mixed module, then $M$ is a endomorphism-liftable if and only if $M=T\oplus F$, where $T$ -- torsion injective module such that all primary components are indecomposable and $F$is a finitely generated projective module.

\textbf{c)} If $M$ is a torsion-free module, then $M$ is endomorphism-liftable if and only if either $M$ is projective or $A$ is a special ring with classical ring of fractions $Q$ and $M=E\oplus F$, where $E$ is a minimal right ideal of the ring $Q$, $F$ is a finitely generated projective module, and $n$ is a positive integer.

\textbf{Proof.} If one of the conditions \textbf{a}, \textbf{b} or \textbf{c} holds, then it follows from Theorem 3.5 that $M$ is an endomorphism-liftable module. In particular, $M$ is an automorphism-liftable, idempotent-liftable module. 

Now let $M$ be an automorphism-liftable, idempotent-liftable module. By Lemma 3.2, $M$ is a $\pi$-projective module. By Theorem 3.4, either one of the conditions \textbf{a} and \textbf{b} of our theorem holds or the following condition holds: there exist two positive integers $k$ and $n$ such that the ring $A$ is isomorphic to the ring of all $k\times k$ matrices $D_k$ over some uniserial Noetherian domain $D$, $M=X\oplus Y$,~ $X$ is a finite direct sum of non-zero injective indecomposable torsion-free modules $X_1, \ldots ,X_n$, $Y$ is a finitely generated projective module and either $n=1$ or $n\ge 2$ and $D$ is a complete domain.

If $n\ge 2$ and $D$ is a complete domain, then $A$ is a special ring, which is required.

Now we assume that $n=1$. Then $A=D$ is a uniserial principal right (left) ideal domain with classical division ring of fractions $Q=E$. Since $Q_A=E_A$ is a direct summand of the automorphism-liftable module $M$, we have that $Q_A$ is an automorphism-liftable module. Since the ring $\End Q_A$ isomorphic to the division ring $Q$, we have that $Q_A$ is an endomorphism-liftable module, by Lemma 3.7. By Theorem 3.5, $A$ is a special ring.~\hfill$\square$

\textbf{3.9. Completion of the proof of Theorem 1.9.} By Lemma 3.3, every automorphism-liftable right $A$-module is an idempotent-liftable module. Therefore, Theorem 1.9 follows from Theorem 3.8.~\hfill$\square$

\section{Remarks and Open Questions} \label{section4}

\textbf{4.1.} There are automorphism-liftable modules which are not endomorphism-liftable; see \cite[Example 5.1]{GQS17} 

\textbf{4.2.} Let $A$ be a non-primitive Dedekind prime ring. If
$\displaystyle{\frac{1}{2}\in A}$, then the automorphism-liftable A-modules coincide with the endomorphism-liftable
A-modules by Theorems 1.9 and 3.5. Are there automorphism-liftable $A$-modules which are not endomorphism-liftable?

\textbf{4.3.} Let $F$ be a field, $A=F[[x]]$ be the formal power series ring, and let $M=F((x))$ be the Laurent series ring. Then  $A$ is a commutative non-primitive Dedekind domain and $M$ is a strongly automorphism-liftable $A$-module which is not singular. 

\textbf{4.4.} For a non-primitive Dedekind prime ring $A$, describe automorphism-liftable $A$-modules and strongly automorphism-liftable $A$-modules. The answer to this question is related to the study invertible elements of the ring $A$ and automorphisms of cyclic $A$-modules.

\textbf{4.5.} By Lemma 3.3, every automorphism-liftable $A$-module is idempotent-liftable provided $A$ contains $\displaystyle{\frac{1}{2}}$. Are there automorphism-liftable modules which are not idempotent-liftable?


\begin{thebibliography}{99}
\bibitem{AbyQ18} Abyzov A.N., Quynh T.C. Lifting of automorphisms of factor modules // Commun. Algebra. - 2018. V.~46, no.~11. -- P.~5073-5082.

\bibitem{GQS17} Guil Asensio P.A., Quynh T.C., Srivastava A.K. // Additive unit structure of endomorphism rings and invariance of modules. -- Bull. Math. Sci. -- 2017. --  Vol.~7. -- P.~229-246.

\bibitem{Jan73} Janakiraman S. Skew projective Abelian groups // Indag. Math. -- 1973. V.76, no.~3. -- P.233--236. 

\bibitem{Len73} Lenagan T.H., Bounded hereditary Noetherian prime rings // J. London Math. Soc. -- 1973. -- Vol.~6. -- P.~241--246. 

\bibitem{MccR87}~McConnell J.~C.,~Robson J.~C. Noncommutative
Noetherian Rings. New York: Wiley-Interscience, 1987.

\bibitem{Mis72} Mishina A. P. On automorphisms and endomorphisms of Abelian groups // Moscow University Mathematics Bulletin. -- 1972. -- no.~1. -- P.~62--66.

\bibitem{SelS18} Selvaraj C., Santhakumar A. S. Automorphism liftable modules // Comment. Math. Univ. Carolin. -- 2018. -- V.~59, no.~1. -- 35-44.

\bibitem{Sin74} Singh S., Quasi-injective and quasi-projective modules over hereditary Noetherian prime rings // Canad. J. Math. -- 1974. -- Vol.~26, no.~5. -- P.~1173--1185.

\bibitem{Sin75}~Singh S. Modules over hereditary Noetherian prime
rings // Can. J. Math. 1975. V.~27, No.~4. P.~867--883.

\bibitem{Tug78} Tuganbaev A.~A. The structure of modules close to projective modules // Sbornik: Mathematics -- 1979. -- V.~35, no.~2. -- P.~219-228.

\bibitem{Tug80} Tuganbaev A.~A. Quasi-projective modules // Sib. math. j. 1980. Vol.~21, no.~3. P.~446--450.

\bibitem{Tug80b} Tuganbaev A.~A. Semiprojective modules // Sib. math. j. 1980. Vol.~21, no.~5. P.~725--728.

\bibitem{Tug01} Tuganbaev A.~A. Modules over bounded Dedekind prime rings // Sbornik: Mathematics -- 2001. -- V.~192, no.~5. -- P.~705-724.

\bibitem{Tug13} Tuganbaev A.~A. Automorphisms of submodules and their extensions // Discrete Math. Appl. -- 2013. -- Vol.~23, no.~1. -- P.~115-124.

\bibitem{Tug16} Tuganbaev A.~A. Automorphism-extendable and
endomorphism-extendable modules // J.~Math. Sci. (New York) -- To appear.

\bibitem{Wis91}~Wisbauer R. Foundations of Module and Ring Theory. Philadelphia: Gordon and Breach, 1991.

\end{thebibliography}
\end{document}